\documentclass{article}
\usepackage{amsmath,amsthm,amssymb,amscd}

\theoremstyle{plain}

\newtheorem{thm}{Theorem}[section]
\theoremstyle{remark}
\newtheorem{rem}{Remark}[section] 
\newtheorem{ex}{Example}[section]
\theoremstyle{definition}

\newcommand{\bb}[1]{\mbox{$\mathbb{#1}$}}

\title{On the dimension formula for the hyperfunction solutions of some holonomic D-modules}
\author{J\"{o}rg Sch\"{u}rmann\thanks{Westf. Wilhelms-Universit\"{a}t, SFB 478 
"Geometrische Strukturen in der Mathematik", Hittorfstr.27, 48149 M\"{u}nster,
Germany,
E-mail: jschuerm@math.uni-muenster.de}}
\date{ }

\begin{document}
\bibliographystyle{plain}

\maketitle

\begin{abstract}
In this short note we improve a recent dimension formula of Takeuchi
for the dimension of the hyperfunction solutions of some holonomic D-modules.
Besides the constructibility result and the local index formula of Kashiwara
for the holomorphic solution complex, we only use a vanishing theorem of Lebeau
together with a simple calculation in terms of constructible functions.\\
$ $
$ $ 
\end{abstract}


\section{Introduction}
One of the basic results about holonomic D-modules is the {\em constructibility result}
of Kashiwara \cite{Ka3}, that the holomorphic solution complex
\[{\cal S}ol({\cal M}):=Rhom_{D_{X}}({\cal M},{\cal O}_{X}) \]
of a holonomic D-module ${\cal M}$ on the complex manifold $X$ is a bounded
{\em complex analytically constructible} complex of sheaves of complex vector spaces with
{\em finite dimensional} stalks (compare also with \cite[chapter 5]{Ka4}, 
\cite[thm.4.5.8, p.458]{KS} and \cite[chap.III]{MNM}).
In particular, the function
\begin{equation} \label{eq:indexdef}
\chi({\cal M}): X\to \bb{Z}\;;\: x\mapsto \chi\bigl(\;{\cal S}ol({\cal M})_{x}\;\bigr) 
\end{equation}
is well defined and {\em complex analytically constructible}. Here $\chi$ is the usual
{\em Euler characteristic}. Moreover, one has by Kashiwara \cite{Ka2} the following beautiful description 
of this local index in terms of the characteristic cycle of the holonomic D-module ${\cal M}$
(see \cite[thm.6.3.1, p.127, cor.6.3.4, p.128]{Ka4} and \cite[thm.2, p.574]{BDK}):

\begin{thm}[{\bf (local index formula)}] \label{thm:index}
Let $Y_{j}$ be finitely many distinct irreducible closed complex analytic subsets of $X$ 
such that the characteristic variety $char({\cal M})$ of ${\cal M}$ can be estimated by
\begin{equation} \label{eq:estimate}
char({\cal M}) \subset \bigcup_{j} \: T^{*}_{Y_{j}}X \:,
\quad \text{with} \quad T^{*}_{Y}X:= cl(T^{*}_{Y_{reg}}X)
\end{equation}
the closure of the conormal bundle to the regular part
$Y_{reg}$ of the irreducible complex analytic subset $Y\subset X$. Then
\begin{equation} \label{eq:index}
\chi({\cal M})(x) = \sum_{j} \: (-1)^{d_{j}} \cdot m_{j} \cdot Eu_{Y_{j}}(x) \:,
\end{equation}
with $d_{j}$ the complex codimension of $Y_{j}$, $Eu_{Y_{j}}$ the famous
Euler obstruction of $Y_{j}$ as defined by MacPherson (cf. \cite{MP}) 
and $m_{j}$ the (generic) multiplicity of ${\cal M}$ along $T^{*}_{Y_{j}}X$.
\end{thm}

Note that $Eu_{Y}=1_{Y}$, if $Y\subset X$ is a {\em closed} complex analytic submanifold.
So a very special case of the local index formula is given as in \cite[ex. on p.129]{Ka4}
by the 

\begin{ex} \label{ex:smooth}
Suppose all $Y_{j}$ in the estimate (\ref{eq:estimate}) are {\em closed} 
connected complex analytic submanifolds of $X$. Then
\begin{equation} \label{eq:indexsmooth}
\chi({\cal M})(x) = \sum_{j} \: (-1)^{d_{j}} \cdot m_{j} \cdot 1_{Y_{j}}(x)\:,
\end{equation}
with $d_{j}$ the complex codimension of $Y_{j}$ and $m_{j}$ the (generic) multiplicity 
of ${\cal M}$ along $T^{*}_{Y_{j}}X$.
\end{ex}

Let us now consider the case that $X$ is the {\em complexification} of the real analytic manifold
$M$, with $i: M\hookrightarrow X$ the closed inclusion. Assume $M$ is purely n-dimensional.
Then the sheaf complex
\[R\Gamma_{M}({\cal O}_{X})[n] \simeq Ri_{*}i^{!}({\cal O}_{X})[n] \]
is concentrated in degree zero, with
\[{\cal B}_{M}:= h^{0}\bigl(\;i^{!}({\cal O}_{X})[n]\;\bigr) \otimes or_{M}\]
the sheaf of Sato's {\em hyperfunctions} on $M$, and $or_{M}$ the {\em orientation sheaf} of $M$.\\

Then the hyperfunction solution complex
\begin{equation} \label{eq:hyperdef}
Rhom_{\,i^{*}D_{X}}(i^{*}{\cal M},{\cal B}_{M}) \simeq 
i^{!}{\cal S}ol({\cal M}) \otimes or_{M}[n]
\end{equation}
of a holonomic D-module ${\cal M}$ on $X$
is {\em subanalytically constructible} on $M$ with {\em finite dimensional} stalks
(compare \cite[thm.5.1.7, p.115]{Ka4}). So it is natural to ask for a corresponding
index formula like (\ref{eq:index}) or (\ref{eq:indexsmooth}).

\begin{rem} The same constructibility result is true for the solutions 
\[Rhom_{\,i^{*}D_{X}}(i^{*}{\cal M},{\cal A}_{M}) \simeq i^{*}{\cal S}ol({\cal M})\]
in the sheaf ${\cal A}_{M}=i^{*}{\cal O}_{X}$
of real analytic functions on $M$.
If ${\cal M}$ is a {\em regular holonomic} D-module on $X$, then one 
also has isomorphisms (see \cite[cor.8.3, cor.8.5, p.360]{Ka4.5} or \cite[p.326]{Bj}):
\[Rhom_{\,i^{*}D_{X}}(i^{*}{\cal M},{\cal B}_{M}) \simeq 
Rhom_{\,i^{*}D_{X}}(i^{*}{\cal M},{\cal D}b_{M})\]
and
\[Rhom_{\,i^{*}D_{X}}(i^{*}{\cal M},{\cal A}_{M}) \simeq 
Rhom_{\,i^{*}D_{X}}(i^{*}{\cal M},{\cal C}^{\infty}_{M}) \:,\]
with ${\cal D}b_{M}$ (or ${\cal C}^{\infty}_{M}$) the sheaf of distributions 
(or smooth functions) on $M$.
\end{rem}

The following counterpart of (\ref{eq:indexsmooth}) is the main result of this note:

\begin{thm}[{\bf (local dimension formula)}] \label{thm:dimension}
Let $M_{j}$ be finitely many distinct closed real analytic submanifolds of $M$ 
such that the characteristic variety $char({\cal M})$ of the holonomic D-module
${\cal M}$ on $X$ can be estimated by
\begin{equation} \label{eq:estimate2}
char({\cal M}) \subset \bigcup_{j} \: T^{*}_{Y_{j}}X \:,
\end{equation}
with $Y_{j}\subset X$ the complexification of $M_{j}$. 
Assume the $Y_{j}$ are irreducible (i.e. connected), with $Y_{j}\cap M=M_{j}$. 
Then one has for $x\in M$:
\begin{equation} \label{eq:dimension}
dim_{\,C}\bigl(\; hom_{\,i^{*}D_{X}}(i^{*}{\cal M},{\cal B}_{M})_{\,x}\;\bigr) 
= \sum_{j} \: m_{j}\cdot 1_{M_{j}}(x)\:,
\end{equation}
with $m_{j}$ the (generic) multiplicity of ${\cal M}$ along $T^{*}_{Y_{j}}X$.
\end{thm}

This is indeed a counterpart of (\ref{eq:indexsmooth}). The estimate (\ref{eq:estimate2})
implies by a theorem of Lebeau \cite{Le} (compare also with \cite[thm.2.1, rem., p.531]{HS}
and \cite[ex.(1), p.533]{HS}) the {\em vanishing result}
\begin{equation} \label{eq:van}
{\cal E}xt^{k}_{\,i^{*}D_{X}}(i^{*}{\cal M},{\cal B}_{M})_{\,x} = 0
\quad \text{for all $k\geq 1$}
\end{equation}
so that
\begin{equation} \label{eq:dim=chi}
dim_{\,C}\bigl(\; hom_{\,i^{*}D_{X}}(i^{*}{\cal M},{\cal B}_{M})_{\,x}\;\bigr) =
\chi\bigl(\;Rhom_{\,i^{*}D_{X}}(i^{*}{\cal M},{\cal B}_{M})_{\,x}\;\bigr) \:.
\end{equation}

Theorem \ref{thm:dimension} answers affirmatively a question asked (or better, discussed)
in \cite[rem.3.5]{Ta} at the end of a recent paper of Takeuchi \cite{Ta},
where he proves the dimension formula (\ref{eq:dimension}) under the  special assumption,
that in suitable local coordinates $(M,x)\simeq (\bb{R}^{n},0)$ the $M_{j}$ are
{\em linear subspaces} (passing through $x=0$).

Note that this special case already covers (locally) the one-dimensional case $(X,M,x)\simeq
(\bb{C},\bb{R},0)$, with ${\cal M}$ a holonomic D-module such that
\[char({\cal M}) \subset T^{*}_{\{x\}}X \cup T^{*}_{X}X \:.\]
In this case one gets back a classical result of Kashiwara \cite[thm.4.2.7, p.69]{Ka1}
(cf. \cite[cor.3.2.36(b), p.88-89]{Ka4}) and Komatsu \cite{Ko}:
\[dim_{\,C}\bigl(\; hom_{\,i^{*}D_{X}}(i^{*}{\cal M},{\cal B}_{M})_{\,x}\;\bigr)
= d + d' \:,\]
with $d$ or $d'$ the multiplicity of ${\cal M}$ along $T^{*}_{\{x\}}X$ or $T^{*}_{X}X$.\\

Let $j: X\backslash M \to M$ be the open inclusion of the complement of $M$ in $X$.
Then the proof given in \cite[sec.3]{Ta} is based on the distinguished triangle
\begin{equation} \begin{CD} \label{eq:triangle}
Ri_{*}i^{!}{\cal S}ol({\cal M})[n] @>>> {\cal S}ol({\cal M})[n] @>>> 
Rj_{*}j^{*}{\cal S}ol({\cal M})[n] @> [1] >> \:.
\end{CD} \end{equation}
Moreover, he uses the {\em micro-local theory} of the {characteristic cycles} $CC(\cdot)$
for subanalytically constructible complexes of sheaves (as in \cite{Ka5, KS, SchVi, Sch}),
in particular a deep result of Schmid-Vilonen \cite{SchVi} about a description of
\[CC(Rj_{*}j^{*}{\cal S}ol({\cal M})) \quad \text{in terms of} \quad
CC({\cal S}ol({\cal M}))\:.\]

In the next section we explain our simple proof of theorem \ref{thm:dimension},
which doesn't make use of this sophisticated micro-local theory of characteristic cycles.
Instead of this, we use the observation that the calculation of
\[\chi(\; (i^{!}{\cal F})_{x}\;) \quad \text{for} \quad {\cal F} = {\cal S}ol({\cal M}) \]
can be done in terms of {\em subanalytically constructible functions}, i.e.
the functor $i^{!}$ induces a corresponding (unique) $\bb{Z}$-linear transformation
for the abelian groups $CF(\cdot)$ of subanalytically constructible functions
such that the following diagram commutes (compare \cite[sec.9.7]{KS} and \cite[sec.2.3]{Sch}):
\begin{equation} \begin{CD} \label{eq:CF}
K_{0}(X) @> i^{!} >> K_{0}(M)\\
@V \chi_{X} VV  @VV \chi_{M} V\\
CF(X) @> i^{!} >> CF(M)\:.
\end{CD} \end{equation}
Here $K_{0}(\cdot)$ is the {\em Grothendieck group} of subanalytically constructible (complexes of)
sheaves with finite dimensional stalks, with $\chi_{?}$ induced by taking stalkwise the Euler
characteristic.\\

Then the calculation of
\begin{equation} \label{eq:i!} 
\chi(\; (i^{!}{\cal S}ol({\cal M})\otimes or_{M}[n])_{x}\;) = 
(-1)^{n}\cdot i^{!}\bigl(\; \chi_{X}({\cal S}ol({\cal M})) \;\bigr)(x)
\end{equation}
becomes an easy exercise by the local index theorem and example \ref{ex:smooth}, since the 
$M_{j}$ and therefore also the $Y_{j}$ are closed submanifolds!\\

If we allow in the estimate (\ref{eq:estimate2}) also {\em singular subspaces},
then we get at least the following weak parity version of the local index theorem:

\begin{thm}[{\bf (local index formula for hyperfunctions)}] \label{thm:indexhyper}
Let $M_{j}$ \\be finitely many distinct real analytic subspaces of $M$ 
such that the characteristic variety $char({\cal M})$ of the holonomic D-module
${\cal M}$ on $X$ can be estimated as in (\ref{eq:estimate2}),
with $Y_{j}\subset X$ the complexification of $M_{j}$. Assume the $Y_{j}$ are
irreducible, with $Y_{j}\cap M=M_{j}$. Then one has for $x\in M$:
\begin{equation} \label{eq:chimod2}
\chi\bigl(\; Rhom_{\,i^{*}D_{X}}(i^{*}{\cal M},{\cal B}_{M})_{\,x}\;\bigr) 
\equiv \sum_{j} \: m_{j}\cdot Eu_{Y_{j}}(x) \quad \text{mod $2$}\:,
\end{equation}
with $Eu_{Y_{j}}$ the Euler obstruction of $Y_{j}$ and
$m_{j}$ the (generic) multiplicity of ${\cal M}$ along $T^{*}_{Y_{j}}X$.
\end{thm}


\section{Constructible functions}

In this final section we give the proof of theorem \ref{thm:dimension} and
\ref{thm:indexhyper} in terms of constructible functions. Let us start with the proof of
theorem \ref{thm:dimension}.\\

By the estimate (\ref{eq:estimate2}) and example \ref{ex:smooth} we get
\[\chi_{X}({\cal S}ol({\cal M})) = \sum_{j} \: (-1)^{d_{j}} \cdot m_{j} \cdot 1_{Y_{j}}\:,\]
with $d_{j}$ the complex codimension of $Y_{j}$ and $m_{j}$ the (generic) multiplicity 
of ${\cal M}$ along $T^{*}_{Y_{j}}X$. By linearity of $i^{!}$ on the level of constructible
functions one also has
\[ i^{!}\bigl(\;\chi_{X}({\cal S}ol({\cal M}))\;\bigr) 
= \sum_{j} \: (-1)^{d_{j}} \cdot m_{j} \cdot i^{!}(1_{Y_{j}}) \:.\]

Then the dimension formula (\ref{eq:dimension}) follows from (\ref{eq:hyperdef}),
(\ref{eq:dim=chi}), (\ref{eq:i!}) and the simple formula
\begin{equation} \label{eq:formula!}
i^{!}(1_{Y_{j}}) = (-1)^{(n-d_{j})}\cdot 1_{M_{j}} \:.
\end{equation}
The formula (\ref{eq:formula!}) corresponds by the commutative diagram (\ref{eq:CF})
to the {\em base change formula}
\[i^{!}(Rk_{*}\bb{C}_{Y_{j}}) \simeq Rk'_{*}i'^{!}\bb{C}_{Y_{j}}\]
for the cartesian diagram of inclusions
\begin{displaymath} \begin{CD}
M @> i >> X \\
@A k' AA @AA k A \\
M_{j} @> i' >> Y_{j} \:.
\end{CD} \end{displaymath}
Note that $i'^{!}\bb{C}_{Y_{j}}\simeq \bb{C}_{M_{j}}[-(n-d_{j})]$ {\em locally} on $M_{j}$,
since $M_{j}$ is a closed submanifold of $Y_{j}$ of real codimension equal to the complex dimension
$n-d_{j}$ of $Y_{j}$.\\

For the proof of the parity formula (\ref{eq:chimod2}) in theorem \ref{thm:indexhyper}
it is enough to show
\begin{equation} \label{eq:mod2}
\chi\bigl(\; i^{*}Rj_{*}j^{*}{\cal S}ol({\cal M})_{x} \;\bigr) 
\equiv 0 \quad \text{mod $2$ for all $x\in M$.}
\end{equation}
Use the local index formula (\ref{eq:index}) and the distinguished triangle
(\ref{eq:triangle}).
But this follows from the fact that the constructible function 
\[\chi_{X}({\cal S}ol({\cal M})) \;mod \;2\] 
is invariant under the complex conjugation acting on 
the complexification $X$ of $M$ (with fixed point set $M$).\\

More precisely, by \cite[lem.1.1.1, p.27]{Sch} one gets the description:
\[\chi\bigl(\; i^{*}Rj_{*}j^{*}{\cal S}ol({\cal M})_{x} \;\bigr)
= \chi\bigl(\;R\Gamma(M_{f,x},{\cal S}ol({\cal M}))\;\bigr)\:,\]
with
\[M_{f,x}:=\{||x||\leq \delta,\; f=w\} \quad \text{for} \quad 0<w<<\delta <<1\]
(i.e. for $w,\delta$ small, with $w$ also small compared to $\delta$) a
{\em local right Milnor fiber} of the function $f$ at $x$, defined in local coordinates
\[(X,M,x) \simeq (\bb{C}^{n},\bb{R}^{n},0) \quad \text{by} \quad
z=(z_{1},\dots,z_{n}) \mapsto f(z):=\sum_{k=1}^{n}\: im(z_{k})^{2} \:.\]
Here $im(\cdot)$ is the imaginary part, with the complex conjugation acting on $(\bb{C}^{n},\bb{R}^{n},0)$
in the usual way. This conjugation leaves the compact semi-analytic set $M_{f,x}$ 
invariant without any fixed point!
But the Euler characteristic
\[\chi\bigl(\;R\Gamma(M_{f,x},{\cal S}ol({\cal M}))\;\bigr) \quad \text{mod $2$}\]
can be calculated in terms of $\bb{Z}_{2}$-valued
constructible functions:
\[\chi\bigl(\;R\Gamma(M_{f,x},{\cal S}ol({\cal M}))\;\bigr) \equiv 
(c\circ \pi)_{*}\alpha \quad \text{mod $2$,}\]
with
\[\alpha := \chi_{M_{f,x}}\bigl(\;{\cal S}ol({\cal M})|M_{f,x}\;\bigr) \quad \text{mod $2$}
\quad \in CF(M_{f,x},\bb{Z}_{2}) \:,\]
\[\pi: M_{f,x}\to M_{f,x}/conj. \quad \text{the quotient and} \quad c: M_{f,x}/conj.\to \{pt\}\] 
a constant map.
Here $(c\circ \pi)_{*}$ is induced by $R(c\circ \pi)_{*}$ similarly as in (\ref{eq:CF})
by the commutative
diagram (compare \cite[sec.9.7]{KS} and \cite[sec.2.3]{Sch}):
\begin{equation} \begin{CD} \label{eq:CF2}
K_{0}(M_{f,x}) @> R(c\circ \pi)_{*} >> K_{0}(\{pt\})\simeq \bb{Z}\\
@V \chi_{M_{f,x}} VV  @VV \chi_{\{pt\}} V\\
CF(M_{f,x}) @> (c\circ \pi)_{*} >> CF(\{pt\})\simeq \bb{Z}\\
@V mod \;2 VV  @VV mod \;2 V\\
CF(M_{f,x},\bb{Z}_{2}) @> (c\circ \pi)_{*} >> CF(\{pt\},\bb{Z}_{2})\simeq \bb{Z}_{2}\:,
\end{CD} \end{equation}  
with $CF(\cdot,\bb{Z}_{2})$ the corresponding 
abelian group of $\bb{Z}_{2}$-valued subanalytically constructible functions.\\

Then $(c\circ \pi)_{*}= c_{*}\circ \pi_{*}$ by functoriality. But
$\pi_{*}(\alpha) \equiv 0$, since $\alpha$ is invariant under the conjugation $conj.$,
with $\pi: M_{f,x}\to M_{f,x}/conj.$ an unramified covering of degree two.
Of course, here it is important to work with $\bb{Z}_{2}$-valued
constructible functions.




\begin{thebibliography}{99}
%
%
\bibitem{Bj}
  Bj\"{o}rk, J.E., \textit{Analytic D-modules and applications},
  Kluwer, 1993.

\bibitem{BDK}
  Brylinski, J.L., Dubson, A. and Kashiwara, M., Formule d'indice pour les
  modules holonomes et obstruction d'Euler locale,
  \textit{C. R. Acad. Sci. Paris} \textbf{293} (1981), 573--576. 

\bibitem{HS}
  Honda, N. and Schapira, P., A vanishing theorem for holonomic modules with
  positive characteristic varieties,
  \textit{Publ. Res. Inst. Math. Sci.} \textbf{26} (1990), 529-534.

\bibitem{Ka1}
  Kashiwara, M., \textit{Algebraic study of systems of linear differential equations}, 
  Master's Thesis, University of Tokyo 1970, English translation in:
  \textit{M\'{e}moires de la SMF} \textbf{63}, 1996.

\bibitem{Ka2}
  Kashiwara, M., Index theorem for maximally overdetermined systems of linear
  differential equations,
  \textit{Proc. Japan Acad.} \textbf{49} (1973), 803--804.

\bibitem{Ka3}
  Kashiwara, M., On the maximally overdetermined systems of linear differential equations~I,
  \textit{Publ. RIMS Kyoto Univ.} \textbf{10} (1975), 563-579.

\bibitem{Ka4}
  Kashiwara, M., \textit{Systems of microdifferential equations},
  \textit{Progress in Mathematics} \textbf{34},
  Birkh\"{a}user, Basel, 1983. 

\bibitem{Ka4.5}
  Kashiwara, M., The Riemann-Hilbert problem for holonomic systems,
  \textit{Publ. RIMS Kyoto Univ.} \textbf{20} (1984), 319-365.

\bibitem{Ka5}
  Kashiwara, M., Index theorem for constructible sheaves,
  \textit{Ast\'{e}risque} \textbf{130} (1985), 193-209.

\bibitem{KS}
  Kashiwara, M. and Schapira, P., \textit{Sheaves on Manifolds},
  Springer, Berlin Heidelberg, 1990.

\bibitem{Ko}
  Komatsu, H., On the index of ordinary differential operators,
  \textit{J. Fac. Sci. Univ. Tokyo Sect. IA Math.} \textbf{18} (1971), 379-398.

\bibitem{Le}
  Lebeau, G., Annulation de la cohomologie hyperfonction de certains modules holonomes,
  \textit{C. R. Acad. Sci. Paris Ser. A-B} \textbf{290} (1980), A313-A316.

\bibitem{MP}
  MacPherson, R., Chern classes for singular varietes,
  \textit{Ann. of Math.} \textbf{100} (1974), 423--432. 

\bibitem{MNM}
  Mebkhout, Z. and Narn\'{a}ez-Macarro, L., Le Th\'{e}or\`{e}me de Constructi\-bilit\'{e}
  de Kashiwara,
  In: \textit{Images directes et constructibilit\'{e}},
  ed. Maisonobe, M. and Sabbah, C., \textit{Travaux en cours},  Hermann, Paris
  \textbf{46} (1997), 47-98.

\bibitem{SchVi}
  Schmid, W. and Vilonen, K., Characteristic cycles for constructible sheaves,
  \textit{Inv. Math.} \textbf{124} (1996), 451-502.

\bibitem{Sch}
  Sch\"{u}rmann, J., \textit{Topology of singular spaces and constructible
sheaves},
  \textit{Monografie Matematyczne (New Series)} \textbf{63}, 
  Birkh\"{a}user, Basel, 2003.

\bibitem{Ta}
  Takeuchi, K., Dimension formulas for the hyperfunction solutions to
  holonomic D-modules,
  \textit{Advances in Mathematics} \textbf{180} (2003), 134-145.


\end{thebibliography}
\end{document}